\documentclass{article}
\usepackage{amssymb}
\usepackage{graphicx}
\usepackage{amsmath}

\input{tcilatex}
\begin{document}

\begin{center}
\ \ \ \ \ \ \ \ \ \ \ \ \ \ \ \ \ \ \ \ \ \ \ \ \ \ \ \ \ \ \ \ \ \ \ \ \ \
\ \ \ \ \ \ \ \ \ \ \ \ \ 

\bigskip

\bigskip \textbf{SINCULAR\ DEGENERATE\ PROBLEMS\ \ AND\ APPLICATIONS}

\ \ \ \textbf{VELI B. SHAKHMUROV}\ \ \ 

Okan University, Department of of Electronics and Communication, Akfirat,
Tuzla 34959 Istanbul, Turkey, E-mail: veli.sahmurov@okan.edu.tr

\ \ \ \ \ \ \ \ \ \ \ \ \ \ \ \ \ \ \ \ \ \ \ 

\textbf{AMS: 34G10,35J25,35J70\ \ }\ \ \ \ \ \ \ \ \ \ \ \ \ 

\textbf{ABSTRACT}
\end{center}

The boundary value problems for linear and nonlinear singular degenerate
differential-operator equations are studied. We prove a well-posedeness of
linear problem and optimal regularity result for the nonlinear problem which
occur in fluid mechanics, environmental engineering and in the atmospheric
dispersion of pollutants.

\begin{center}
\textbf{Key Words: }differential-operator equations, Semigroups of
operators, Banach-valued function spaces, separability, fredholmness,
interpolation of Banach spaces, Atmospheric dispersion of pollutants.

\ \ \textbf{1. Introduction, notations and background }
\end{center}

The maximal regularity properties for boundary value problems (BVPs) for
linear differential-operator equations (DOEs) have been studied extensively
by many researchers (see e.g. $\left[ 1-10\right] $ and the references
therein). The main objective of the present paper is to discuss BVPs for the
following nonlinear singular degenerate DOE 
\begin{equation*}
-x^{2\alpha }\frac{\partial ^{2}u}{\partial x^{2}}-y^{2\beta }\frac{\partial
^{2}u}{\partial y^{2}}+a\left( x,u,u_{x},u_{y}\right) u=F\left(
x,u,u_{x},u_{y}\right) 
\end{equation*}%
on the rectangular domain $G=\left( 0,a\right) \times \left( 0,b\right) .$

Several conditions for the uniform separability and the resolvent estimates
for the corresponding linear problem are given in abstract $L_{p}$-spaces.
Especially, we prove that the linear differential operator is positive and
is a generator of an analytic semigroup. Moreover, the existence and
uniqueness of maximal regular solution of the above nonlinear problem are
obtained. One of the important characteristics of these DOEs are that the
degeneration process are taking place at different speeds at boundary, in
general. Maximal regularity properties of regular degenerated nonlinear DOEs
are studied e.g. in $\left[ \text{1, 8, 10}\right] .$ Unlike to these we
consider here the singular degenerate DOEs. In applications maximal
regularity properties of infinite systems of singular degenerate PDE are
studied.

Let $\gamma =\gamma \left( x\right) ,$ $x=\left(
x_{1},x_{2},...,x_{n}\right) $ be a positive measurable function on a domain 
$\Omega \subset R^{n}.$ Let $L_{p,\gamma }\left( \Omega ;E\right) $ denote
the space of strongly measurable $E$-valued functions that are defined on $%
\Omega $ with the norm

\begin{equation*}
\left\Vert f\right\Vert _{L_{p,\gamma }}=\left\Vert f\right\Vert
_{L_{p,\gamma }\left( \Omega ;E\right) }=\left( \int \left\Vert f\left(
x\right) \right\Vert _{E}^{p}\gamma \left( x\right) dx\right) ^{\frac{1}{p}},%
\text{ }1\leq p<\infty .
\end{equation*}

For $\gamma \left( x\right) \equiv 1$, $L_{p,\gamma }\left( \Omega ;E\right) 
$ will be denoted by $L_{p}=L_{p}\left( \Omega ;E\right) $. Let $\mathbf{C}$
be the set of the complex numbers and\ 
\begin{equation*}
S_{\varphi }=\left\{ \lambda ;\text{ \ }\lambda \in \mathbf{C}\text{, }%
\left\vert \arg \lambda \right\vert \leq \varphi \right\} \cup \left\{
0\right\} ,\text{ }0\leq \varphi <\pi .
\end{equation*}

A linear operator\ $A$ is said to be $\varphi $-positive in a Banach\ space $%
E$ with bound $M>0$ if $D\left( A\right) $ is dense on $E$ and $\left\Vert
\left( A+\lambda I\right) ^{-1}\right\Vert _{B\left( E\right) }\leq M\left(
1+\left\vert \lambda \right\vert \right) ^{-1}$ $\ $for any $\lambda \in
S_{\varphi },$ $0\leq \varphi <\pi ,$ where $I$ is the identity operator in $%
E$ and $B\left( E\right) $ is the space of bounded linear operators in $E.$

The $\varphi $-positive operator $A$ is said to be $R$-positive in a Banach
space $E$ if the set $L_{A}=\left\{ \xi \left( A+\xi I\right) ^{-1}\text{: }%
\xi \in S_{\varphi }\right\} ,$ $0\leq \varphi <\pi $ is $R$-bounded (see
e.g. $\left[ 5\right] $ ). Let $G$ be a domain in $R^{n}$. Let\ $W_{p,\gamma
}^{m}=W_{p,\gamma }^{m}\left( G;E\left( A\right) ,E\right) $ and $%
W_{p,\gamma }^{\left[ m\right] }=W_{p,\gamma }^{\left[ m\right] }\left(
0,a;E\left( A\right) ,E\right) $ are $E$-valued weighted function spaces
defined in $\left[ 8\right] .$

\begin{center}
\textbf{2. Linear degenerate DOEs}
\end{center}

\ Consider the BVP for the singular degenerate differential-operator
equation 
\begin{equation}
-x^{2\alpha }\frac{\partial ^{2}u}{\partial x^{2}}-x^{2\beta }\frac{\partial
^{2}u}{\partial y^{2}}+Au+x^{\alpha }A_{1}\frac{\partial u}{\partial x}%
+y^{\beta }A_{2}\frac{\partial u}{\partial y}+\lambda u=f\left( x,y\right) ,
\end{equation}

\begin{equation*}
L_{1}u=\sum\limits_{i=0}^{m_{1}}\delta _{1i}u_{x}^{\left[ i\right] }\left(
a,y\right) =0,L_{2}u=\sum\limits_{i=0}^{m_{2}}\delta _{2i}u_{y}^{\left[ i%
\right] }\left( x,b\right) =0,\text{ }
\end{equation*}%
on the domain $G=\left( 0,a\right) \times \left( 0,b\right) ,$ where $%
u=u\left( x,y\right) ,$ $u_{x}^{\left[ i\right] }=\left[ x^{\alpha }\frac{%
\partial }{\partial x}\right] ^{i}u\left( x,y\right) ,$ $u_{y}^{\left[ i%
\right] }=\left[ y^{\beta }\frac{\partial }{\partial y}\right] ^{i}u\left(
x,y\right) $, $m_{k}\in \left\{ 0,1\right\} $; $\delta _{jik}$ are complex
numbers, $\lambda $ is a complex parameter, $A$ and $A_{i}=A_{i}\left(
x,y\right) $ are linear operators in a Banach space $E.$

Let $\delta _{1m_{k}}\neq 0,$ $k=1,2.$The main result is the following

\textbf{Theorem 1. }Let $E$ be an $UMD$ space space (see$\left[ 11\right] $
), $A$ be an $R$-positive operator in $E,$ $A_{i}A^{-\left( \frac{1}{2}-\mu
\right) }\in L_{\infty }\left( G;B\left( E\right) \right) $ for $\mu \in
\left( 0,\frac{1}{2}\right) $ and $1+\frac{1}{p}<\alpha ,\beta <\frac{\left(
p-1\right) }{2}$. Then the problem $\left( 1\right) $ has a unique solution $%
u\in W_{p,\alpha ,\beta }^{\left[ 2\right] }\left( G;E\left( A\right)
,E\right) $ for all  $f\in L_{p}\left( G;E\right) $ and $\left\vert \arg
\lambda \right\vert \leq \varphi $ with sufficently large $\left\vert
\lambda \right\vert $ and the following coercive uniform estimate holds

\begin{equation}
\sum\limits_{i=0}^{2}\left\vert \lambda \right\vert ^{1-\frac{i}{2}}\left[
\left\Vert x^{i\alpha }\frac{\partial ^{i}u}{\partial x^{i}}\right\Vert
_{L_{p}}+\left\Vert y^{i\beta }\frac{\partial ^{i}u}{\partial y^{i}}%
\right\Vert _{L_{p}}\right] +\left\Vert Au\right\Vert _{L_{p}}\leq
M\left\Vert f\right\Vert _{L_{p}}.
\end{equation}

For proving the main theorem, consider at fist, BVPs for the singular
degenerate DOE

\begin{equation}
\ \left( L+\lambda \right) u=-u^{\left[ 2\right] }\left( x\right) +\left(
A+\lambda \right) u\left( x\right) =f,
\end{equation}%
\begin{equation*}
L_{1}u=\sum\limits_{i=0}^{m_{k}}\delta _{i}u^{\left[ i\right] }\left(
a\right) =0\text{, }
\end{equation*}%
where $u^{\left[ i\right] }=\left[ x^{\alpha }\frac{d}{dx}\right]
^{i}u\left( x\right) $, $m_{k}\in \left\{ 0,1\right\} ;$ $\delta _{i}$ are
complex numbers and $A$ is a linear operator in $E,$ $\delta _{m_{k}}\neq 0$.

In a similar way as in $\left[ \text{9, Theorem 5.1}\right] $ we obtain

\textbf{Theorem A}$_{1}$\textbf{. }Suppose\textbf{\ }$E$ is an $UMD$ space, $%
A$ is an $R$ positive in $E$, $1+\frac{1}{p}<\alpha <\frac{\left( p-1\right) 
}{2}.$ Then the problem $\left( 3\right) $ has a unique solution $u\in
W_{p,\alpha }^{\left[ 2\right] }$ for all \ $f\in L_{p}\left( 0,a;E\right) ,$
$p\in \left( 1,\infty \right) .$ Moreover for $\left\vert \arg \lambda
\right\vert \leq \varphi $ and sufficiently large $\left\vert \lambda
\right\vert $ the following uniform coercive estimate holds

\begin{equation}
\sum\limits_{i=0}^{2}\left\vert \lambda \right\vert ^{1-\frac{i}{2}%
}\left\Vert u^{\left[ i\right] }\right\Vert _{L_{p}\left( 0,a;E\right)
}+\left\Vert Au\right\Vert _{L_{p}\left( 0,a;E\right) }\leq C\left\Vert
f\right\Vert _{L_{p}\left( 0,a;E\right) }.
\end{equation}

Let $B$ denote the operator in $L_{p}\left( 0,a;E\right) $ generated by
problem $\left( 2\right) $, i.e. 
\begin{equation*}
D\left( B\right) =\left\{ u:u\in W_{p,\alpha }^{\left[ 2\right] },\text{ }%
L_{k}u=0\right\} ,\text{ }Bu=-u^{\left[ 2\right] }+Au.
\end{equation*}
In a similar way as in $\left[ \text{8, Theorem 3.1}\right] $ we obtain

\textbf{Theorem A}$_{2}.$ Let all conditions of Theorem A$_{1}$ are
satisfied. Then, the operator $B$ is $R$-positive in $L_{p}\left(
0,a;E\right) .$

Theorem A$_{1}$ implies that the operator $B$ is positive and is a generator
of analytic semigroups in $L_{p}\left( 0,a;E\right) $.

Consider now the following degenerate DOEs with the boundary conditions $%
\left( 3\right) :$%
\begin{equation}
\ -x^{2\alpha }u^{\left( 2\right) }\left( x\right) +\left( A+\lambda \right)
u\left( x\right) =f,\text{ }L_{k}u=0.
\end{equation}

\textbf{Theorem A}$_{3}.$ Let all conditions of Theorem A$_{1}$ are
satisfied. Then the problem $\left( 5\right) $ has a unique solution $u\in
W_{p,\alpha }^{2}$ for all \ $f\in L_{p}\left( 0,a;E\right) .$ Moreover for $%
\left\vert \arg \lambda \right\vert \leq \varphi $ and sufficiently large $%
\left\vert \lambda \right\vert $ the following coercive estimate holds

\begin{equation}
\sum\limits_{i=0}^{2}\left\vert \lambda \right\vert ^{1-\frac{i}{2}%
}\left\Vert x^{i\alpha }u^{\left( i\right) }\right\Vert _{L_{p}\left(
0,a;E\right) }+\left\Vert Au\right\Vert _{L_{p}\left( 0,a;E\right) }\leq
C\left\Vert f\right\Vert _{L_{p}\left( 0,a;E\right) }.
\end{equation}

\textbf{Proof. }Since $\alpha >1,$ by $\left[ \text{9, Theorem 2.3}\right] $
we get that there is a small $\varepsilon >0$ and $C\left( \varepsilon
\right) $ such that 
\begin{equation}
\left\Vert \alpha x^{\alpha -1}u^{\left[ 1\right] }\right\Vert _{L_{p}\left(
0,a;E\right) }\leq \varepsilon \left\Vert u\right\Vert _{W_{p,\alpha }^{%
\left[ 2\right] }\left( 0,a;E\left( A\right) ,E\right) }+C\left( \varepsilon
\right) \left\Vert u\right\Vert _{L_{p}\left( 0,a;E\right) }.
\end{equation}%
Then in view of $\left( 6\right) $, $\left( 7\right) $ and due to positivity
of operator $B,$ we have the following estimate 
\begin{equation}
\left\Vert \alpha x^{\alpha -1}u^{\left[ 1\right] }\right\Vert _{L_{p}\left(
0,a;E\right) }\leq \varepsilon \left\Vert Bu\right\Vert _{L_{p}\left(
0,a;E\right) }.
\end{equation}

Since $-x^{2\alpha }u^{\left( 2\right) }=-u^{\left[ 2\right] }+\alpha
x^{\alpha -1}u^{\left[ 1\right] },$ the assertion is obtained from Theorem A$%
_{1}$ and the estimate $\left( 8\right) .$

In this stage we can show the proof of Theorem1.

\textbf{Proof of Theorem 1.} Consider at first the principal part of the
problem $\left( 1\right) $ i.e

\begin{equation}
-x^{2\alpha }\frac{\partial ^{2}u}{\partial x^{2}}-x^{2\beta }\frac{\partial
^{2}u}{\partial y^{2}}+Au+\lambda u=f\left( x,y\right) ,\text{ }L_{1}u=0%
\text{, }L_{2}u=0.
\end{equation}

Since $L_{p}\left( 0,b;L_{p}\left( 0,a;E\right) \right) =$ $L_{p}\left(
G;E\right) $ then, this BVP can be express as:

\begin{equation}
-y^{2\beta }\frac{d^{2}u}{dy^{2}}+\left( B+\lambda \right) u\left( y\right)
=f\left( y\right) \text{, }L_{2k}u=0.
\end{equation}

By virtue of $\left[ \text{1, Theorem 4.5.2}\right] $ $F=L_{p}\left(
0,b;E\right) \in UMD$ provided $E\in UMD$, $p\in \left( 1,\infty \right) $.
By Theorem A$_{2},$\ the operator $B\ $is $R$-positive in $F.$ Then by
virtue of Theorem A$_{3}$, for $f\in L_{p}\left( 0,a;F\right) =L_{p}\left(
G;E\right) $ problem $\left( 9\right) $ has a unique solution $u\in $ $%
W_{p,\beta }^{2}\left( 0,a;D\left( S\right) ,F\right) $ and \ the operator $Q
$ generated by problem $\left( 9\right) $ has a bounded inverse from $%
L_{p}\left( G;E\right) $ to $W_{p,\alpha ,\beta }^{\left[ 2\right] }.$
Moreover by using embedding theorems in $W_{p,\alpha ,\beta }^{\left[ 2%
\right] }$ (see e.g. $\left[ \text{9, Theorem 2.3}\right] $ we get the
following estimate%
\begin{equation*}
\left\Vert x^{\alpha }A_{1}\frac{\partial u}{\partial x}\right\Vert
_{L_{p}\left( G;E\right) }+\left\Vert y^{\beta }A_{2}\frac{\partial u}{%
\partial y}\right\Vert _{L_{p}\left( G;E\right) }\leq \varepsilon \left\Vert
Qu\right\Vert _{L_{p}\left( G;E\right) }\text{, }\varepsilon <1.
\end{equation*}%
By virtue the above estimate and by using perturbation properties of linear
operators we obtain the assertion.

\textbf{Remark 1. }Note that, by using the similar techniques similar to
those applied in Theorems 1, 2, we can obtained the same results for
differential-operator equations of the arbitrary order.

\begin{center}
3. \textbf{Singular degenerate BVPs with small parameters}
\end{center}

\bigskip\ Consider the BVP for the parameter dependent degenerate
differential-operator equation 

\begin{equation}
\ Lu=-tu^{\left[ 2\right] }\left( x\right) +\left( A+\lambda \right) u\left(
x\right) =f,\text{ }x\in \left( 0,a\right) 
\end{equation}%
\begin{equation*}
L_{1}u=\sum\limits_{i=0}^{1}t^{\sigma _{i}}\alpha _{i}u^{\left[ i\right]
}\left( a\right) =f_{1},
\end{equation*}%
where $u^{\left[ i\right] }=\left[ x^{\gamma }\frac{d}{dx}\right]
^{i}u\left( x\right) ,$ $\sigma _{i}=\frac{i}{2}+\frac{1}{2\left( 1-\gamma
\right) p},$ $\gamma >1+\frac{1}{p}$; $\alpha _{i}$ are complex numbers $t$
is a small positive and $\lambda $ is a complex parameter, $A$ is a linear
operator in a Banach space $E$ and  $f_{1}\in E_{1}=\left( \left( E\left(
A\right) ,E\right) _{\theta }\right) ,$ $\theta =\frac{1}{2}\left( 1+\frac{1%
}{\left( 1-\alpha \right) p}\right) .$ \ 

A function \ $u\in $ $W_{p,\gamma }^{\left[ 2\right] }\left( 0,a;E\left(
A\right) ,E\right) $\ satisfying the equation $\left( 1\right) $ a.e. on $%
\left( 0,1\right) $ is said to be the solution of the equation $\left(
1\right) $ on $\left( 0,1\right) .$

\textbf{Remark 2.}\ Let 
\begin{equation}
y=\int\limits_{0}^{x}z^{-\gamma }dz.
\end{equation}

Under the substitution $\left( 11\right) $ spaces $L_{p}\left( 0,1;E\right) $
and $W_{p,\gamma }^{\left[ 2\right] }\left( 0,a;E\left( A\right) ,E\right) $
are mapped isomorphically onto weighted spaces $L_{p,\tilde{\gamma}}\left(
-\infty ,0;E\right) $ and 
\begin{equation*}
W_{p,\tilde{\gamma}}^{2}\left( -\infty ,0;E\left( A\right) ,E\right) ,\text{ 
}\tilde{\gamma}=\tilde{\gamma}\left( x\left( y\right) \right) .
\end{equation*}%
Moreover, under the substitution $\left( 10\right) $ the problem $\left(
2\right) $ is transformed into a non degenerate problem

\begin{equation}
Lu=-tu^{\left( 2\right) }\left( y\right) +Au\left( y\right) =f,\text{ }%
L_{1}u=\sum\limits_{i=0}^{1}t^{\sigma _{i}}\alpha _{i}u^{\left( i\right)
}\left( 0\right) =f_{1},
\end{equation}

\textbf{Theorem 2. }Suppose\ $E$ is a UMD and $1+\frac{1}{p}<\gamma
,1<p<\infty .$ Then the problem $\left( 10\right) $ for all \ $f\in
L_{p}\left( 0,a;E\right) ,$ $f_{1}\in E_{1}$ has a unique solution $u\in
W_{p,\gamma }^{\left[ 2\right] }\left( 0,a;E\left( A\right) ,E\right) $ and
for $\left\vert \arg \lambda \right\vert \leq \varphi $ and sufficiently
large $\left\vert \lambda \right\vert $ the following uniform coercive
estimate holds

\begin{equation}
\sum\limits_{i=0}^{2}\left\vert \lambda \right\vert ^{1-\frac{i}{2}}t^{\frac{%
i}{2}}\left\Vert u^{\left[ i\right] }\right\Vert _{L_{p}\left( 0,a;E\right)
}+\left\Vert Au\right\Vert _{L_{p}\left( 0,a;E\right) }\leq C\left[
\left\Vert f\right\Vert _{L_{p}\left( 0,a;E\right) }+\left\Vert
f_{1}\right\Vert _{E_{1}}+\left\vert \lambda \right\vert ^{1-\theta
}\left\Vert f_{1}\right\Vert _{E}.\right] 
\end{equation}

\textbf{Proof. }\ Consider the problem $\left( 12\right) $. In a similar way
as in $\left[ \text{10, Theorem 3.2}\right] $ we obtain that the problem $%
\left( 12\right) $ has a unique solution $u\in W_{p,\tilde{\gamma}}^{\left(
2\right) }\left( -\infty ,0;E\left( A\right) ,E\right) $ for all \ $f\in
L_{p,\tilde{\gamma}}\left( -\infty ,0;E\right) ,$ $f_{1}\in E_{1}$ and $%
\left\vert \arg \lambda \right\vert \leq \varphi $ with sufficiently large $%
\left\vert \lambda \right\vert $ the following uniform coercive estimate
holds

\begin{equation*}
\sum\limits_{i=0}^{2}\left\vert \lambda \right\vert ^{1-\frac{i}{2}}t^{\frac{%
i}{2}}\left\Vert u^{\left( i\right) }\right\Vert _{L_{p,\tilde{\gamma}%
}\left( -\infty ,0;E\right) }+\left\Vert Au\right\Vert _{L_{p,\tilde{\gamma}%
}\left( -\infty ,0;E\right) }\leq C\left[ \left\Vert f\right\Vert _{L_{p,%
\tilde{\gamma}}\left( -\infty ,0;E\right) }+\left\Vert f_{1}\right\Vert
_{E_{1}}+\left\vert \lambda \right\vert ^{1-\theta }\left\Vert
f_{1}\right\Vert _{E}.\right] 
\end{equation*}

Then in view of the Remark 2 we obtain the assertion.

Consider now the parameter dependent singular degenerate BVP%
\begin{equation}
-t_{1}x^{2\alpha }\frac{\partial ^{2}u}{\partial x^{2}}-t_{2}x^{2\beta }%
\frac{\partial ^{2}u}{\partial y^{2}}+Au+\lambda u=f\left( x,y\right) ,
\end{equation}

\begin{equation*}
L_{1}u=\sum\limits_{i=0}^{m_{1}}t_{1}^{\sigma _{1i}}\delta _{1i}u_{x}^{\left[
i\right] }\left( a,y\right) =0,L_{2}u=\sum\limits_{i=0}^{m_{2}}t_{2}^{\sigma
_{2i}}\delta _{2i}u_{y}^{\left[ i\right] }\left( x,b\right) =0,\text{ }
\end{equation*}%
on the domain $G=\left( 0,a\right) \times \left( 0,b\right) $, where $\sigma
_{1i}=\frac{i}{2}+\frac{1}{2p\left( 1-\alpha \right) },$ $\sigma _{2i}=\frac{%
i}{2}+\frac{1}{2p\left( 1-\beta \right) }$,  $A$ is a linear operator in a
Banach space $E$ and $t_{k}$ are small parameters.

\textbf{Theorem 3. }Let $E$ be an $UMD$ space space, $A$ be an $R$-positive
operator in and $1+\frac{1}{p}<\alpha ,\beta <\frac{\left( p-1\right) }{2}$.
Then the problem $\left( 14\right) $ has a unique solution $u\in W_{p,\alpha
,\beta }^{\left[ 2\right] }\left( G;E\left( A\right) ,E\right) $ for all $%
f\in L_{p}\left( G;E\right) $ and $\left\vert \arg \lambda \right\vert \leq
\varphi $, with sufficiently large $\left\vert \lambda \right\vert $ and the
following coercive uniform estimate holds

\begin{equation}
\sum\limits_{i=0}^{2}\left\vert \lambda \right\vert ^{1-\frac{i}{2}}\left[
t_{1}^{i}\left\Vert x^{i\alpha }\frac{\partial ^{i}u}{\partial x^{i}}%
\right\Vert _{L_{p}}+t_{2}^{i}\left\Vert y^{i\beta }\frac{\partial ^{i}u}{%
\partial y^{i}}\right\Vert _{L_{p}}\right] +\left\Vert Au\right\Vert
_{L_{p}}\leq M\left\Vert f\right\Vert _{L_{p}}.
\end{equation}

\textbf{Proof. }\ By reasoning as in the proof of Theorem1, the problem $%
\left( 14\right) $ is reduced to \ the following BVP for ordinary equation 
\begin{equation}
\ -t_{2}u^{\left[ 2\right] }\left( y\right) +\left( B_{t_{1}}+\lambda
\right) u\left( y\right) =f,\text{ }y\in \left( 0,b\right) 
\end{equation}%
\begin{equation*}
L_{2}u=\sum\limits_{i=0}^{1}t_{2}^{\sigma _{i}}\alpha _{i}u^{\left[ i\right]
}\left( b\right) =0,
\end{equation*}

where $B_{t_{2}}$ is the operator in $L_{p}\left( 0,b;E\right) $\ generated
by BVP 
\begin{equation*}
\ Lu=-t_{1}u^{\left[ 2\right] }\left( x\right) +\left( A+\lambda \right)
u\left( x\right) =f,\text{ }x\in \left( 0,a\right) 
\end{equation*}%
\begin{equation*}
L_{1}u=\sum\limits_{i=0}^{1}t_{1}^{\sigma _{i}}\alpha _{i}u^{\left[ i\right]
}\left( a\right) =0.
\end{equation*}

Then by applying Theorem 2 to problem $\left( 16\right) $ in $L_{p}\left(
0,a;F\right) =L_{p}\left( G;E\right) ,$ $F=L_{p}\left( 0,b;E\right) $ we
obtain the assertion.

\begin{center}
4. \textbf{Singular degenerate BVPs in moving domains}
\end{center}

Consider the linear BVPs in moving domain $G_{s}=\left( 0,a\left( s\right)
\right) \left( 0,b\left( s\right) \right) $%
\begin{equation}
-x^{2\alpha }\frac{\partial ^{2}u}{\partial x^{2}}-x^{2\beta }\frac{\partial
^{2}u}{\partial y^{2}}+Au+du=f\left( x,y\right) ,
\end{equation}

\begin{equation*}
L_{1}u=\sum\limits_{i=0}^{1}\delta _{1i}u_{x}^{\left[ i\right] }\left(
a\left( s\right) ,y\right) =0,L_{2}u=\sum\limits_{i=0}^{1}\delta _{2i}u_{y}^{%
\left[ i\right] }\left( x,b\left( s\right) \right) =0,\text{ }k=1,2,
\end{equation*}%
\qquad \qquad \qquad where $a\left( s\right) $ and $b(s)$ are positive
continues function depended on parameter $s$ and $A$ is a linear operator in
a Banach space $E,$ $\delta _{k1}$ $\neq 0.$\ 

\textbf{Theorem 4. }Let $E$ be an $UMD$ space space, $A$ be an $R$-positive
operator in $E,$ $1+\frac{1}{p}<\alpha ,\beta <\frac{\left( p-1\right) }{2}$
and $f_{k}\in E_{k}$.

Then problem $\left( 10\right) $\ for $f\in L_{p}\left( G\left( s\right)
;E\right) $, $f_{k}\in E_{k},$ $p\in \left( 1,\infty \right) $ and the
sufficiently large $d>0$ has a unique solution $u\in $\ $W_{p,\alpha ,\beta
}^{2}\left( G\left( s\right) ;E\left( A\right) ,E\right) $ and the following
coercive uniform estimate holds

\begin{equation*}
\left\Vert x^{2\alpha }\frac{\partial ^{2}u}{\partial x^{2}}\right\Vert
_{L_{p}\left( G\left( s\right) ;E\right) }+\left\Vert y^{2\beta }\frac{%
\partial ^{2}u}{\partial y^{2}}\right\Vert _{L_{p}\left( G\left( s\right)
;E\right) }+\left\Vert Au\right\Vert _{L_{p}\left( G\left( s\right)
;E\right) }\leq 
\end{equation*}

\begin{equation*}
C\left\Vert f\right\Vert _{L_{p}\left( G_{s};E\right) }.
\end{equation*}%
\textbf{Proof. }Under the substitution $\tau =xb(s),$ $\sigma =yb\left(
s\right) $ by denoting $\tau $, $\sigma ,$ $u\left( x\left( \tau \right)
,y\left( \sigma \right) \right) ,$ $f\left( x\left( \tau \right) ,y\left(
\sigma \right) \right) $ again by $x,$ $y,$\ $u\left( x,y\right) ,$ $f\left(
x,y\right) ,$ respectively we get the moving boundary problem $\left(
10\right) $ maps to the following BVP with parameter in the fixed domain $%
\left( 0,1\right) \times \left( 0,1\right) $%
\begin{equation}
-b^{2\left( 1-\alpha \right) }\left( s\right) x^{2\alpha }\frac{\partial
^{2}u}{\partial x^{2}}-b^{2\left( 1-\beta \right) }\left( s\right) x^{2\beta
}\frac{\partial ^{2}u}{\partial y^{2}}+Au+du=f\left( x,y\right) ,
\end{equation}

\begin{equation*}
L_{1}u=\sum\limits_{i=0}^{1}b^{i}\left( s\right) \delta _{1i}u_{x}^{\left[ i%
\right] }\left( 1,y\right) =0,L_{2}u=\sum\limits_{i=0}^{1}b^{i}\left(
s\right) \delta _{2i}u_{y}^{\left[ i\right] }\left( x,1\right) =0,\text{ }%
k=1,2.
\end{equation*}

Then by virtue of Theorem 3 we obtain the assertion.

\begin{center}
\textbf{5}. \textbf{Nonlinear degenerate DOE}
\end{center}

Consider now the following nonlinear problem%
\begin{equation}
-x^{2\alpha }\frac{\partial ^{2}u}{\partial x^{2}}-y^{2\beta }\frac{\partial
^{2}u}{\partial y^{2}}+a\left( u,u_{x},u_{y}\right) u=F\left(
x,y,u,u_{x},u_{y}\right) ,\text{ }L_{1k}u=0,\text{ }L_{2k}u=0\text{ }
\end{equation}%
on the domain $G_{0}=\left( 0,a_{0}\right) \times \left( 0,b_{0}\right) ,$
where $L_{jk}$ are boundary conditions defined by $\left( 1\right) .$ Let 
\begin{equation*}
X_{1}=L_{p}\left( 0,a;E\right) \text{, }Y_{1}=W_{p,\alpha }^{\left[ 2\right]
}\left( 0,a;E\left( A\right) ,E\right) ,\text{ }X_{2}=L_{p}\left(
0,b;E\right) ,
\end{equation*}%
\begin{equation*}
Y_{2}=W_{p,\beta }^{\left[ 2\right] }\left( 0,b;E\left( A\right) ,E\right) ,%
\text{ }E_{ki}=\left( X_{k},Y_{k}\right) _{\theta _{ki},p},\text{ }\theta
_{1i}=\frac{p\left( 1-\alpha \right) i+1}{2p\left( 1-\alpha \right) },\text{ 
}
\end{equation*}

\begin{equation*}
\theta _{2i}=\frac{p\left( 1-\beta \right) i+1}{2p\left( 1-\beta \right) },%
\text{ }E_{0}=\prod\limits_{i,k}E_{ki},\text{ }i=0,1,\text{ }k=1,2.
\end{equation*}%
\textbf{Condition 1. }Assume the following satisfied:

(1)\ $E$ is an $UMD$ space$,$ $a\left( x,y,U\right) =A\left( x,y\right) $ is
a positive operator in $E$ for $x,y\in G_{0}$, $u_{i}\in E_{1i},$ $g_{i}\in
E_{1i}$, $D\left( a\left( x,y,U\right) \right) $ does not depend on $x,y,U$,
where $U=\left\{ u_{0},u_{1},g_{0},g_{1}\right\} $ and $a:$ $G_{0}\times
E_{0}\rightarrow B\left( E\left( A\right) ,E\right) $ is continuous;

(2) $F:G_{0}\times E_{0}\rightarrow E$ be a measurable function; $F\left(
x,y,.\right) $ is continuous with respect to $x,y\in G_{0}$ and $f\left(
x,y\right) =$ $F\left( x,y,0\right) \in X.$ Moreover, for each $R>0$ there
exists $\mu _{R}$ such that $\left\Vert F\left( x,U\right) -F\left( x,\bar{U}%
\right) \right\Vert _{E}\leq \mu _{R}\left\Vert U-\bar{U}\right\Vert
_{E_{0}} $\ for a.a. $x,y\in G_{0}$, $u_{j},$ $\bar{u}_{j}\in X_{j}$ and $%
\left\Vert U\right\Vert _{E_{0}}\leq R,\left\Vert \bar{U}\right\Vert
_{E_{0}}\leq R$, $1+\frac{1}{p}<\alpha ,\beta <\frac{\left( p-1\right) }{2}$;

(3) there exist $v_{j}\in E_{1j}$, $\upsilon _{j}\in E_{2j}$ such that the
operator $A\left( x,y,\Phi \right) $ for $\Phi =\left\{ v_{1},v_{2},\upsilon
_{1},\upsilon _{2}\right\} $ is $R$-positive in $E$ uniformly with respect
to $x,y\in G_{0};$ $A\left( x,y,\Phi \right) A^{-1}\left( x^{0},y^{0},\Phi
\right) \in C\left( G_{0};B\left( E\right) \right) ;$

(4) Moreover, for each $R>0$ there is a positive constant $L\left( R\right) $
such that

$\left\Vert \left[ A\left( x,y,U\right) -A\left( x,y,\bar{U}\right) \right]
\upsilon \right\Vert _{E}\leq M\left( R\right) \left\Vert U-\bar{U}%
\right\Vert _{E_{0}}\left\Vert A\upsilon \right\Vert _{E}$ for $x,y\in G_{0}$%
, $\left\Vert U\right\Vert _{E_{0}},\left\Vert \bar{U}\right\Vert
_{E_{0}}\leq R$ and $\upsilon \in D\left( A\left( x,y,U\right) \right) .$

In this section we prove the existence and uniqueness of maximal regular
solution for the nonlinear problem $\left( 11\right) $. 

\textbf{Theorem 3. }Let the Condition1 holds. Then there is $a\in \left(
0\right. \left. a_{0}\right] ,$ $b\in \left( 0\right. \left. b_{0}\right] $
such that problem $\left( 11\right) $ has a unique solution belongs to $%
W_{p,\alpha ,\beta }^{2}\left( (G;E\left( A\right) ,E\right) .$

\textbf{Proof. }By Theorem , the linear problem 
\begin{equation}
\ -x^{2\alpha }\frac{\partial ^{2}w}{\partial x^{2}}-y^{2\beta }\frac{%
\partial ^{2}w}{\partial y^{2}}+Aw\left( x,y\right) =f\left( x,y\right) ,
\end{equation}%
\begin{equation*}
L_{1k}u=0,\text{ }L_{2k}u=0,\text{ }k=1,2,\text{ }x,y\in \left( 0,a\right)
\times \left( 0,b\right) 
\end{equation*}%
is maximal regular in $X$ \ uniformly with respect to $a\in \left( 0\right.
\left. a_{0}\right] $ and $b\in \left( 0\right. \left. b_{0}\right] $ i.e.
for all $f\in X$ there is a unique solution $w\in Y$ of the problem $(31)$
and has a coercive estimate 
\begin{equation*}
\left\Vert w\right\Vert _{Y}\leq C\left\Vert f\right\Vert _{X},
\end{equation*}%
where the constant $C$\ does not depends on $a\in \left( 0\right. \left.
a_{0}\right] $ and 
\begin{equation*}
f\left( x\right) =F\left( x,0\right) .
\end{equation*}%
We want to to solve the problem $\left( 30\right) $ locally by means of
maximal regularity of the linear problem $(31)$ via the contraction mapping
theorem. For this purpose let $w$ be a solution of the linear BVP $(31).$
Consider a ball 
\begin{equation*}
B_{r}=\left\{ \upsilon \in Y,\text{ }\left\Vert \upsilon -w\right\Vert
_{Y}\leq r\right\} .
\end{equation*}

Given $\upsilon \in B_{r},$ solve the problem%
\begin{equation}
-tu^{\left( 2m\right) }\left( x\right) +Au\left( x\right) =F\left(
x,\upsilon ,\upsilon ^{\left( 1\right) },...,\upsilon ^{\left( 2m-1\right)
}\right) ,\text{ }
\end{equation}%
\begin{equation*}
\sum\limits_{i=0}^{m_{k}}t^{\eta _{i}}\left[ \alpha _{ki}u^{\left( i\right)
}\left( 0\right) +\beta _{ki}u^{\left( i\right) }\left( a\right)
+\sum\limits_{j=1}^{N_{k}}\delta _{kj}u^{\left( i\right) }\left(
x_{kj}\right) \right] =f_{k},\text{ }k=1,2,....,2m
\end{equation*}%
where $x\in \left( 0,a\right) .$ Define a map $Q$ on $B_{r}$ by $Q\upsilon
=u,$ where $u$ is a solution of the problem $\left( 38\right) .$ We want to
show that $Q\left( B_{r}\right) \subset B_{r}$ and that $L$ is a contraction
operator in $Y$, provided $a$ is sufficiently small, and $r$ is chosen
properly. For this aim by using maximal regularity properties of the problem 
$\left( 37\right) $ for $V=\left\{ \upsilon ^{\left( m_{k}\right) }\left(
0\right) \right\} ,$ $k=1,2,...,2m$ we have 
\begin{equation*}
\left\Vert Q\upsilon -w\right\Vert _{Y}=\left\Vert u-w\right\Vert _{Y}\leq
C_{0}\left\Vert F\left( x,V\right) -F\left( x,0\right) \right\Vert _{X}.
\end{equation*}

By assumption Condition1 and in view of Remark1 we have 
\begin{equation*}
\left\Vert F\left( x,V\right) -F\left( x,0\right) \right\Vert _{E}\leq
\end{equation*}%
\begin{equation*}
\left\Vert F\left( x,V\right) -F\left( x,W\right) \right\Vert
_{E}+\left\Vert F\left( x,W\right) -F\left( x,0\right) \right\Vert _{E}\leq
\end{equation*}

\begin{equation*}
M_{R}\left[ \left\Vert V-W\right\Vert _{E_{0}}+\left\Vert W\right\Vert
_{E_{0}}\right] 
\end{equation*}%
\begin{equation*}
M_{R}C_{1}\left[ \left\Vert \upsilon -w\right\Vert _{Y}+\left\Vert
w\right\Vert _{Y}\right] \leq MC_{1}\left[ r+\left\Vert w\right\Vert _{Y}%
\right] ,
\end{equation*}%
where $R=M_{R}C_{1}\left[ r+\left\Vert w\right\Vert _{Y}\right] $ is a fixed
number such that $R\leq \frac{r}{C_{0}}$. In view of Condition C and the
above estimates for sufficiently small $a\in \lbrack 0;a_{0})$we have 
\begin{equation*}
\left\Vert Q\upsilon -w\right\Vert _{Y}\leq C_{0}R\leq r
\end{equation*}%
i.e. 
\begin{equation*}
Q\left( B_{r}\right) \subset B_{r}.
\end{equation*}%
In a similar way for $\bar{V}=\left\{ \bar{\upsilon}^{\left( m_{k}\right)
}\left( 0\right) \right\} $ we obtain 
\begin{equation*}
\left\Vert Q\upsilon -Q\bar{\upsilon}\right\Vert _{Y}\leq C_{0}\left[
\left\Vert F\left( x,V\right) -F\left( x,\bar{V}\right) \right\Vert _{X}%
\right] \leq C_{0}M_{R}\left\Vert \left( \upsilon -\bar{\upsilon}\right)
\right\Vert _{Y}.
\end{equation*}%
Therefore for $C_{0}M_{R}<1$ the operator $Q$ becomes a contraction mapping.
Eventually, the contraction mapping principle implies a unique fixed point
of $Q$ in $B_{r}$ which is the unique strong solution $u\in
Y=W_{p}^{2m}\left( 0,a;E\left( A\right) ,E\right) .$ 

\bigskip 

\begin{center}
\textbf{\ 4. Singular degenerate boundary value problems for infinite
systems of equations }
\end{center}

Consider the infinite system of BVPs%
\begin{equation}
-x^{2\alpha }\frac{\partial ^{2}u_{m}}{\partial x^{2}}-x^{2\beta }\frac{%
\partial ^{2}u_{m}}{\partial y^{2}}+d_{m}u_{m}+\dsum\limits_{j=1}^{\infty
}x^{\alpha }a_{mj}\left( x,y\right) \frac{\partial u_{j}}{\partial x}
\end{equation}%
\begin{equation*}
+\dsum\limits_{j=1}^{\infty }y^{\beta }b_{mj}\left( x,y\right) \frac{%
\partial u_{j}}{\partial y}+\lambda u=f_{m}\left( x,y\right) ,\text{ }%
L_{1k}u=0,\text{ }L_{2k}u=0,
\end{equation*}

where $L_{ik}$ are defined by $\left( 1\right) $. Let $D$ 
\begin{equation*}
D=\left\{ d_{m}\right\} ,\text{ }d_{m}>0,\text{ }u=\left\{ u_{m}\right\} ,%
\text{ }Du=\left\{ d_{m}u_{m}\right\} ,\text{ }m=1,2,...,
\end{equation*}

\begin{equation*}
\text{ }l_{q}\left( D\right) =\left\{ u\text{: }u\in l_{q},\right.
=\left\Vert u\right\Vert _{l_{q}\left( D\right) }=\left. \left(
\sum\limits_{m=1}^{\infty }\left\vert d_{m}u_{m}\right\vert ^{q}\right) ^{%
\frac{1}{q}}<\infty ,q\in \left( 1,\infty \right) \right\} .
\end{equation*}
From Theorem1we obtain

\textbf{Theorem 3. }Assume $a_{mj},b_{mj}\in L_{\infty }\left( G\right) $.
For $0<\mu <\frac{1}{2}$ and for all $x,y\in \left( G\right) $ 
\begin{equation*}
\text{ }\sup\limits_{m}\sum\limits_{j=1}^{\infty }a_{mj}\left( x\right)
d_{j}^{-\left( \frac{1}{2}-\mu \right) }<M,\text{ }\sup\limits_{m}\sum%
\limits_{j=1}^{\infty }b_{mj}\left( x\right) d_{j}^{-\left( \frac{1}{2}-\mu
\right) }.
\end{equation*}%
Then for all $f\left( x\right) =\left\{ f_{m}\left( x\right) \right\}
_{1}^{\infty }\in L_{p}\left( \left( G\right) ;l_{q}\right) ,$ $p,q\in
\left( 1,\infty \right) $, $\left\vert \arg \lambda \right\vert \leq \varphi 
$, $0\leq \varphi <\pi $ and for sufficiently large $\left\vert \lambda
\right\vert $ problem $\left( 12\right) $ has a unique solution $u=\left\{
u_{m}\left( x\right) \right\} _{1}^{\infty }$ that belongs to space $%
W_{p,\alpha ,\beta }^{2}\left( G,l_{q}\left( D\right) ,l_{q}\right) $ and 
\begin{eqnarray}
&&\sum\limits_{i=0}^{2}\left\vert \lambda \right\vert ^{1-\frac{i}{2}}\left[
\left\Vert x^{i\alpha }\frac{\partial ^{i}u}{\partial x^{i}}\right\Vert
_{L_{p}\left( G;l_{q}\right) }+\left\Vert y^{i\beta }\frac{\partial ^{i}u}{%
\partial y^{i}}\right\Vert _{L_{p}\left( G;l_{q}\right) }\right] +\left\Vert
Du\right\Vert _{L_{p}\left( G;l_{q}\right) } \\
&\leq &M\left\Vert f\right\Vert _{L_{p}\left( G;l_{q}\right) }.  \notag
\end{eqnarray}

\textbf{References}\ \ \ \ \ \ \ \ \ \ \ \ \ \ \ \ \ \ \ \ \ \ \ \ \ \ \ \ \
\ \ \ \ \ \ \ \ \ \ \ \ \ \ \ \ \ \ \ \ \ \ \ \ \ \ \ \ \ \ \ \ \ \ \ \ \ \
\ \ \ \ \ \ \ \ \ \ \ \ \ \ \ \ \ 

\begin{enumerate}
\item Amann H., Linear and quasi-linear equations,1, Birkhauser, Basel 1995.

\item Yakubov S. and Yakubov Ya., \textquotedblright Differential-operator
Equations. Ordinary and Partial \ Differential Equations \textquotedblright
, Chapman and Hall /CRC, Boca Raton, 2000.

\item Krein S. G., \textquotedblright Linear differential equations in
Banach space\textquotedblright , American Mathematical Society, Providence,
1971.

\item Sobolevskii P. E., Coerciveness inequalities for abstract parabolic
equations, Doklady Akademii Nauk SSSR, 57(1),(1964), 27-40.

\item Denk R., Hieber M., Pr\"{u}ss J., $R$-boundedness, Fourier multipliers
and problems of elliptic and parabolic type, Mem. Amer. Math. Soc. 166
(2003), n.788.

\item Favini A., Shakhmurov V., Yakubov Y., Regular boundary value problems
for complete second order elliptic differential-operator equations in UMD
Banach spaces, Semigroup form, v. 79 (1), 2009, 22-54.

\item Ashyralyev A., Claudio Cuevas and Piskarev S., "On well-posedness of
difference schemes for abstract elliptic problems in spaces", Numerical
Functional Analysis \& Optimization, v. 29, No. 1-2, Jan. 2008, 43-65.

\item Shakhmurov V. B, Shahmurova A., Nonlinear abstract boundary value
problems atmospheric dispersion of pollutants, Nonlinear Analysis, Wold
Applications, v.11 (2) 2010, 932-951.

\item Shakhmurov V. B., Degenerate differential operators with parameters,
Abstract and Applied Analysis, 2007, v. 2007, 1-27.

\item Shakhmurov V. B., Nonlinear abstract boundary value problems in
vector-valued function spaces and applications, Nonlinear Analysis Series A:
Theory, Method \& Applications, v. 73, 2010, 2383-2397.

\item Burkholder D. L., A geometrical conditions that implies the existence
certain singular integral of Banach space-valued functions, Proc. conf.
Harmonic analysis in honor of Antonu Zigmund, Chicago, 1981,Wads Worth,
Belmont, (1983), 270-286.

\item Triebel H., \textquotedblright Interpolation theory, Function spaces,
Differential operators.\textquotedblright , North-Holland, Amsterdam, 1978.
\end{enumerate}

\end{document}